\documentclass[10pt]{amsart}
\usepackage{amsfonts,amssymb}
\def\qq{\mathbb{Q}}
\def\rr{\mathbb{R}}
\def\cc{\mathbb{C}}
\def\zz{\mathbb{Z}}

\def\Ardeg{\widehat{\deg}}

\def\difftials{H^{0}(X,\Omega_{X}^{1})}
\def\imtau{\mathrm{Im}\tau}

\def\ee{\mathcal{E}}

\def\isom{{\buildrel \sim \over \longrightarrow}}

\theoremstyle{plain}
\newtheorem{theorem}{Theorem}[section]
\newtheorem{lemma}[theorem]{Lemma}
\newtheorem{proposition}[theorem]{Proposition}
\newtheorem{corollary}[theorem]{Corollary}

\theoremstyle{remark}
\newtheorem{remark}[theorem]{Remark}

\theoremstyle{definition}
\newtheorem{definition}[theorem]{Definition}

\theoremstyle{definition}
\newtheorem{definition/proposition}[theorem]{Definition/Proposition}

\numberwithin{equation}{section}

\begin{document}

\title{On the Arakelov theory of elliptic curves}
\author{Robin de Jong}
\address{University of Amsterdam, The Netherlands}
\email{rdejong@science.uva.nl}

\begin{abstract}
This note contains an elementary discussion of the Arakelov
intersection theory of elliptic curves. The main new results are a
projection formula for elliptic arithmetic surfaces and a formula
for the ``energy'' of an isogeny between Riemann surfaces of genus 1.
The latter formula provides an answer to a question originally
posed by Szpiro.
\end{abstract}

\date{\today}

\maketitle

\section{Introduction}

The goal of this note is to present an
elementary discussion of the Arakelov intersection theory of
elliptic curves. Arakelov intersection theory in general is a
theory dealing with curves over number fields, unifying in a
subtle way the arithmetic aspects of a curve, present on the
reductions of the curve modulo the finite primes of the number
field, with its analytic aspects, present on the Riemann surfaces
that one obtains by base changing the curve to the complex
numbers. The unifying framework is provided by an intersection
theory for divisors on an arithmetic surface \cite{ar}, sharing
many formal properties with the traditional intersection theory on
proper algebraic surfaces over a field \cite{fa}. Although in
general working out Arakelov theory is a difficult matter, when
we specify to the case of elliptic curves it turns out that a
nice, compact and clean theory emerges.

Many results on the Arakelov theory of elliptic curves 
are already known by the works of Faltings \cite{fa} and 
Szpiro \cite{sz}, but our approach is different. In particular, we base
our discussion on a projection formula for Arakelov's Green
function on Riemann surfaces of genus 1 related by an isogeny. From
this formula we derive a projection formula for Arakelov intersections, as well as a formula for the so-called
``energy of an isogeny''. Both of these formulas seem new. In
fact, the latter formula provides an answer to a question posed by
Szpiro in \cite{sz}.

Using these new results, we give alternative proofs of several of
the earlier results. For example, we arrive at explicit formulas
for the Arakelov-Green function on an elliptic curve, for the
canonical norm in the holomorphic cotangent bundle, and for the
Arakelov self-intersection of a point. We also give an
elementary proof of a recent result due to Autissier on the
average height of the quotients of an elliptic curve by its cyclic
subgroups of a fixed order.

\section{Analytic invariants} \label{preliminaries}

We start by recalling the main ingredients of the analytic part of
Arakelov theory, namely, the Arakelov-Green function $G$ and the
canonical metric on the holomorphic cotangent bundle. Our main
references are \cite{ar} and \cite{fa}.

Let $X$ be a compact and connected Riemann surface of genus $g>0$.
The space of holomorphic differentials $\difftials$ carries then a
natural hermitian inner product $(\omega,\eta) \mapsto \frac{i}{2}
\int_X \omega \wedge \overline{\eta}$. Let $\{
\omega_1,\ldots,\omega_g \}$ be an orthonormal basis with respect
to this inner product. We have then a fundamental (1,1)-form $\mu$
on $X$ given by $ \mu = \frac{i}{2g} \sum_{k=1}^g \omega_k \wedge
\overline{\omega}_k $. It is verified immediately that the form
$\mu$ does not depend on the choice of orthonormal basis, and
hence it defines a canonical (1,1)-form on $X$. Using this form,
one defines the canonical Arakelov-Green function on $X$. This
function gives the local intersections ``at infinity'' of two
divisors in Arakelov theory.
\begin{definition}
The Arakelov-Green function $G$ is
the unique function $ X \times X \to \rr_{\geq 0}$ such that
the following three properties hold:
\begin{itemize} \item[(i)] $G(P,Q)^2$ is $C^\infty$ on $X \times X$
and $G(P,Q)$ vanishes only at the diagonal $\Delta_X$, with multiplicity
1;
\item[(ii)] for all $P \in X$ we have $\partial_Q
\overline{\partial}_Q
\log G(P,Q)^2 = 2\pi i \mu(Q)$ for $Q \neq P$;
\item[(iii)] for all $P \in X$ we have $\int_X \log G(P,Q) \mu(Q)=0$.
\end{itemize}
Properties (i) and (ii) determine $G$ up to a multiplicative
constant, which is then fixed by the normalisation condition
(iii). By an application of Stokes' theorem we obtain from
(i)--(iii) the symmetry $G(P,Q)=G(Q,P)$ of the function $G$.
\end{definition}

Importantly, the Arakelov-Green function gives rise to certain
canonical metrics on the line bundles $O_X(D)$, where $D$ is a
divisor on $X$. It suffices to consider the case of a point $P \in
X$, for the general case follows then by taking tensor products.
Let $s$ be the canonical generating section of the line bundle
$O_X(P)$. We then define a smooth hermitian metric $\| \cdot
\|_{O_X(P)}$ on $O_X(P)$ by putting $\|s\|_{O_X(P)}(Q) := G(P,Q)$
for any $Q \in X$. By property (ii) of the Arakelov-Green
function, the curvature form of $O_X(P)$ is equal to $\mu$, and in
general, the curvature form of $O_X(D)$ is $\deg(D) \cdot \mu$,
with $\deg(D)$ the degree of $D$.
\begin{definition} A line bundle $L$ on $X$ with a smooth hermitian
metric $\| \cdot \|$ is called admissible if its curvature form is
a multiple of $\mu$. We also call the metric $\| \cdot \|$ itself
admissible in this case.
\end{definition}
\begin{proposition} \label{quotconstant}
Let $\| \cdot \|$ and $\| \cdot \|'$ be
admissible metrics on a line bundle $L$. Then the quotient $\|
\cdot \|/\| \cdot \|'$ is a constant function on $X$.
\end{proposition}
\begin{proof} The logarithm of the quotient is a smooth harmonic
function on $X$, and hence it is constant.
\end{proof}
\begin{definition}
The canonical metric $\| \cdot \|_{\mathrm{Ar}}$ on the holomorphic
cotangent bundle $\Omega_X^1$ is the unique metric that makes the 
adjunction isomorphism $O_{X \times X}(-\Delta_X)|_{\Delta_X}
\isom \Omega_X^1$ an isometry. Here the line bundle $O_{X\times
X}(\Delta_X)$ carries the hermitian metric defined by
$\|s\|(P,Q):=G(P,Q)$, with $s$ the canonical generating section of the
line bundle $O_{X \times X}(\Delta_X)$.
\end{definition}
\begin{proposition} \label{adjunction} (Adjunction formula) Let $P$ be a
point on $X$,
and let $z$ be a local coordinate about $P$. Then for the norm
$\|dz\|_{\mathrm{Ar}}$ of $dz$ in $\Omega_X^1$ the formula $\| dz
\|_{\mathrm{Ar}} = \lim_{Q \to P} |z(P)-z(Q)|/G(P,Q)$ holds.
\end{proposition}
\begin{proof} From the definition of the canonical metric on
$\Omega_X^1$ it follows that $dz/z$ has unit length in
$\Omega_X^1(P)$. However, this line bundle is isometric to $
\Omega_X^1 \otimes O_X(P)$, with $dz/z$ corresponding to $dz
\otimes   z^{-1}s$ with $s$ the canonical generating section of
$O_X(P)$. One computes that $\|   z^{-1}s \| = \lim_{Q \to P}
G(P,Q)/|z(P)-z(Q)|$ and the proposition follows.
\end{proof}
In Sections \ref{compl} and \ref{kernelofisogeny} we prove some
fundamental properties of the Arakelov-Green function and the
canonical norm on the holomorphic cotangent bundle in the case
that $X$ has genus~1.

\section{Analytic projection formula} \label{compl}

We start by studying the fundamental (1,1)-form $\mu$ with respect to
isogenies. Let $X$ and $X'$ be
Riemann surfaces of genus 1, and suppose that $f : X \to X'$ is an
isogeny, say of degree $N$. Let $\mu_X$ and $\mu_{X'}$ be the
fundamental (1,1)-forms of $X$ and $X'$, respectively.
\begin{proposition} \label{prelim}
(i) We have $f^* \mu_{X'} = N \cdot \mu_X$;
(ii) the canonical isomorphism $f^* : H^0(X',\Omega^1_{X'}) \isom
\difftials$ given by inclusion has norm $\sqrt{N}$.
\end{proposition}
\begin{proof} We identify $X$ with a complex torus $\cc/\Lambda$,
and obtain $X'$ as the quotient of $\cc/\Lambda$ by a finite
subgroup $\Lambda'/\Lambda$. Hence we may identify $X'$ with
$\cc/\Lambda'$. A small computation shows that the differentials
$\omega := dz/\sqrt{\mathrm{vol}(\Lambda)}$ and $\omega' :=
dz/\sqrt{\mathrm{vol}(\Lambda')}$ are orthonormal bases of
$\difftials$ and $H^0(X',\Omega^1_{X'})$, respectively. We obtain
(ii) by observing that $N =
\mathrm{vol}(\Lambda)/\mathrm{vol}(\Lambda')$. Finally we have
$\mu_X = (i/2) \cdot (dz \wedge
d\overline{z})/\mathrm{vol}(\Lambda)$ and $\mu_{X'} = (i/2) \cdot
(dz \wedge d\overline{z})/\mathrm{vol}(\Lambda')$ and (i) also
follows.
\end{proof}
Proposition \ref{prelim} gives rise to a projection formula for
the Arakelov-Green function.
\begin{theorem} \label{complexproj} (Analytic projection formula)
Let $X$ and $X'$ be Riemann surfaces of genus
1 and let $G_X$ and $G_{X'}$ be the Arakelov-Green functions of
$X$ and $X'$, respectively. Suppose we have an isogeny $f : X \to
X'$. Let $D$ be a divisor on $X'$. Then the canonical isomorphism
of line bundles
\[ f^* O_{X'}(D) \isom O_X(f^*D) \] is an isometry. In particular
we have a projection formula: for any $P \in X$ the formula
\[ G_X(f^*D,P) = G_{X'}(D,f(P)) \] holds.
\end{theorem}
\begin{proof} Let $N$ be the degree of $f$. By Proposition
\ref{prelim} we have
\[ \begin{aligned}
\mathrm{curv} f^*O_{X'}(D) & = f^*( \mathrm{curv} O_{X'}(D)) =
f^*( (\deg D) \cdot \mu_{X'}) = N \cdot (\deg D) \cdot \mu_X \\ &=
\deg(O_X(f^*D)) \cdot \mu_X \, ,
\end{aligned}
\] which means that $f^*O_{X'}(D)$ is an admissible line bundle on
$X$. Hence by Proposition \ref{quotconstant} we have $ \| f^*(s_D)
\|_{f^*O_{X'}(D)} = c \cdot \| s_{f^*D} \|_{O_X(f^*D)}$ for some
constant $c$ where $s_D$ and $s_{f^*D}$ are the canonical sections
of $O_{X'}(D)$ and $O_X(f^*D)$, respectively. But since
\[ \begin{aligned}
\int_X \log &\| f^*(s_D) \|_{f^*O_{X'}(D)} \cdot \mu_X =
\frac{1}{N} \cdot \int_X \log \| f^*(s_D) \|_{f^*O_{X'}(D)} \cdot
f^*\mu_{X'}  \\ &= \int_{X'} \log \| s_D \|_{O_{X'}(D)} \cdot
\mu_{X'} = 0 \, , \end{aligned} \] this constant is equal to 1.
\end{proof}

\section{Energy of an isogeny} \label{kernelofisogeny}

At this point, we introduce some classical invariants attached to
a Riemann surface $X$ of genus 1.
\begin{definition} Let $\tau$ be an element
of the complex upper half plane, and write $q=\exp(2\pi i \tau)$.
Then we have the eta-function $\eta(\tau) = q^{1/24}
\prod_{k=1}^\infty (1-q^k)$ and the modular discriminant
$\Delta(\tau) = \eta(\tau)^{24} = q \prod_{k=1}^\infty
(1-q^k)^{24}$. The latter is the unique normalised cusp form of
weight 12 on $\mathrm{SL}(2,\zz)$. Now suppose that we have a
Riemann surface $X$ of genus 1 identified with a complex torus
$\cc/\zz+\tau \zz$. Then we put $\|\eta\|(X) := (\imtau)^{1/4}
\cdot |\eta(\tau)|$ and $\|\Delta\|(X)
:=\|\eta\|(X)^{24}=(\imtau)^6 \cdot |\Delta(\tau)|$. These
definitions do not depend on the choice of $\tau$, and hence they
define invariants of $X$.
\end{definition}

In \cite{sz} Szpiro proves 
the following statement (\emph{cf.}
Th\'eor\`eme 1): let $E$ and $E'$ be semi-stable elliptic curves
defined over a number field $K$, and suppose we have an isogeny $f
: E \to E'$. Then the formula
\[ \sum_\sigma \sum_{{P_\sigma \in \mathrm{Ker} f_\sigma, \atop P_\sigma
\neq 0}} \log G(0,P_\sigma) = \frac{ [K:\qq] }{2} \log N +
\sum_\sigma \log \frac{ \| \eta \|(E'_\sigma)^2 }{ \| \eta
\|(E_\sigma)^2 } \] holds, where $N$ is the degree of $f$ and
where the sum is over the complex embeddings of $K$. 
Szpiro then asks whether a similar statement holds without the
sum over the complex embeddings. The following theorem gives a
positive answer to that question. The terminology ``energy of an
isogeny'' is adopted from \cite{sz}.
\begin{theorem} \label{kernel} (Energy of an isogeny)
Let $X$ and $X'$ be Riemann surfaces of genus 1.
Suppose we have an isogeny $f : X \to X'$. Then we have
\[ \prod_{P \in \mathrm{Ker} f, P \neq 0} G(0,P) = \frac{ \sqrt{N}
\cdot \| \eta \|(X')^2 }{ \| \eta \|(X)^2 } \, , \] where $N$ is
the degree of $f$.
\end{theorem}
It is the purpose of the present section to prove Theorem
\ref{kernel}. \emph{En passant} we make the Arakelov-Green
function and the canonical norm on the holomorphic cotangent
bundle explicit, see Propositions \ref{explG} and \ref{AX}. These
formulas are also given in \cite{fa}, but the proof there relies
on a consideration of the eigenvalues and eigenfunctions of the
Laplace operator. Our approach is more elementary.
\begin{definition} Let $X$ be a Riemann surface of genus 1. Let
$\omega$ be a holomorphic differential of norm 1 in $\difftials$.
Then we put $A(X) := \|\omega\|_{\mathrm{Ar}}$ for the norm of
$\omega$ in $\Omega_X^1$.
\end{definition}
\begin{proposition} \label{roughkernel}
Let $f : X \to X'$ be an isogeny of degree $N$. Then the formula
\[ \prod_{P \in \mathrm{Ker} f, P \neq 0} G(0,P) = \frac{ \sqrt{N}
\cdot A(X) }{A(X')} \] holds.
\end{proposition}
\begin{proof} Let $\nu$ be the norm of the isomorphism of line bundles
$f^* \Omega^1_{X'} {\buildrel \sim \over \to} \Omega^1_X$ given by
the usual inclusion. We will compute $\nu$ in two ways. First of
all, consider an $\omega' \in H^0(X',\Omega^1_{X'})$ of norm 1, so
that $\omega'$ has norm $A(X')$ in $\Omega^1_{X'}$. Then by
Proposition \ref{prelim} we have that $f^*(\omega')$ has norm
$\sqrt{N}$ in $\difftials$, hence it has norm $\sqrt{N} \cdot
A(X)$ in $\Omega_X^1$. This gives
\[ \nu = \frac{ \sqrt{N} \cdot A(X) }{A(X')} \, . \]
On the other hand, by Theorem \ref{complexproj}, the canonical
isomorphism $ f^*(O_{X'}(0)) {\buildrel \sim \over \to}
O_X(\mathrm{Ker} f)$ is an isometry. Tensoring with the
isomorphism $f^* \Omega^1_{X'} {\buildrel \sim \over \to}
\Omega^1_X $ gives an isomorphism
\[ f^*(\Omega_{X'}^1(0)) {\buildrel \sim \over \longrightarrow} \Omega_X^1(0) \otimes
\bigotimes_{P \in \mathrm{Ker} f, P \neq 0} O_X(P) \] of norm
$\nu$  given in local coordinates by
\[ f^*( \frac{dz}{z}) \mapsto \frac{dz}{z} \otimes s \]
where $s$ is the canonical section of $\bigotimes_{P \in
\mathrm{Ker} f, P \neq 0} O_X(P) $. By the definition of the
canonical norm on the holomorphic cotangent bundle, the $dz/z$
have norm 1, so we find
\[ \nu = \prod_{P \in \mathrm{Ker} f, P \neq 0} G(0,P) \, . \]
Together with the earlier formula for $\nu$ this implies the
proposition.
\end{proof}
The following corollary seems to be well-known, see for instance
\cite{szpull}, Lemme 6.2.
\begin{corollary} \label{Ntorsion}
Denote by $X[N]$ the kernel of the
multiplication-by-$N$ map $X \to X$. Then the formula
\[ \prod_{P \in X[N], P \neq 0} G(0,P) = N \] holds.
\end{corollary}
\begin{proof} Immediate from Proposition \ref{roughkernel}.
\end{proof}
Let $\tau$ be an element of the complex upper half plane. Recall
that Riemann's theta function is given by $\vartheta(z;\tau):=
\sum_{n \in \zz} \exp(\pi i n^2 \tau +2\pi i nz)$ on $\cc$. We have the
identities
\[ \left( \exp(\pi i \tau/4) \cdot \vartheta(0;\tau)
\vartheta(1/2;\tau) \vartheta(\tau/2;\tau) \right)^8 = 2^8 \cdot
\Delta(\tau) \] and
\[ \left( \exp(\pi i \tau/4) \cdot \frac{ \partial \vartheta}{\partial z}
\left( \frac{1+\tau}{2};\tau \right) \right)^8 = (2\pi)^8 \cdot
\Delta(\tau) \, , \] both of 
which are proved by the fact that the left hand
sides are cusp forms on $\mathrm{SL}(2,\zz)$ of weight 12.
\begin{definition} (\emph{Cf.} \cite{fa}) Let $\tau$ be in the
complex upper half plane. The normalised theta function $\|
\vartheta \|$ associated to $\tau$ is defined to be the function
\[ \|\vartheta\|(z;\tau) := (\imtau)^{1/4} \exp(-\pi (\imtau)^{-1} y^2)
|\vartheta(z;\tau)| \] on $\cc$ where $y := \mathrm{Im} z$. This
function only depends on the class of $z$ modulo $\zz+\tau\zz$.
\end{definition}
\begin{proposition} \label{explG} (Faltings \cite{fa})
Let $X$ be a Riemann surface of genus 1, and
write $X \cong \cc/\zz+\tau \zz$ with $\tau$
in the complex upper half plane. For the Arakelov-Green function $G$ on $X$ the
formula
\[ G(0,z) = \frac{ \| \vartheta \| (z + (1+\tau)/2;\tau) }{ \| \eta
\|(X) }   \] holds.
\end{proposition}
\begin{proof} It is not difficult to check that
$\|\vartheta\|(z+(1+\tau)/2)$ vanishes only at $z=0$, 
with order 1. Also it is not difficult to check
that $\partial_z \overline{\partial}_z \log
\|\vartheta\|(z+(1+\tau)/2)^2 = 2\pi i \mu_X$ for $z \neq 0$. By
what we have said in Section \ref{preliminaries}, we have from
this that $G(0,z)=c \cdot \| \vartheta \| (z + (1+\tau)/2;\tau)$
where $c$ is some constant. It remains to compute this constant.
If we apply Corollary \ref{Ntorsion} with $N=2$ we obtain
\[ c^3 \cdot \| \vartheta \| (0;\tau)
\| \vartheta \| (1/2;\tau)  \| \vartheta \| (\tau/2;\tau) =
G(0,1/2)G(0,\tau/2)G(0,(1+\tau)/2) = 2 \, . \] On the other hand
we have the formula
\[ \left( \exp(\pi i \tau/4) \cdot \vartheta(0;\tau)
\vartheta(1/2;\tau) \vartheta(\tau/2;\tau) \right)^8 = 2^8 \cdot
\Delta(\tau) \,  \] mentioned above. Combining we obtain
$ c = \|\eta\|(X)^{-1}$.
\end{proof}
\begin{proposition} \label{AX} (Faltings \cite{fa})
For the invariant $A(X)$, the
formula
\[ A(X) = \frac{1}{ (2\pi) \cdot \| \eta \|(X)^2 }   \] holds.
\end{proposition}
\begin{proof} We follow the argument from \cite{fa}: writing $X
\cong \cc/\zz+\tau \zz$ we can take $\omega = dz/\sqrt{\imtau}$ as
an orthonormal basis of $\difftials$. By Proposition
\ref{adjunction} we have $ \| dz/\sqrt{\imtau} \|_{\mathrm{Ar}} =
(\sqrt{\imtau})^{-1} \cdot \lim_{z \to 0} |z| / G(0,z)$. We obtain
the required formula by using the explicit formula for $G(0,z)$ in
Proposition \ref{explG} and the formula
\[ \left( \exp(\pi i \tau/4) \cdot \frac{ \partial \vartheta}{\partial z}
\left( \frac{1+\tau}{2};\tau \right) \right)^8 = (2\pi)^8 \cdot
\Delta(\tau)  \] mentioned above.
\end{proof}
\begin{proof}[Proof of Theorem \ref{kernel}] 
Immediate from Propositions
\ref{roughkernel} and \ref{AX}.
\end{proof}
We conclude this section with a corollary, dealing with the value
of the Arakelov-Green function on pairs of 2-torsion points. First
we need a classical lemma.
\begin{lemma} \label{classicfactsI} Let $X$ be a Riemann surface of genus 1 and suppose
that $y^2=4x^3-px-q=: f(x)$ is a Weierstrass equation for $X$.
Write $f(x) =4(x-\alpha_1)(x-\alpha_2)(x-\alpha_3)$. Let
$(\omega_1|\omega_2)$ be the period matrix of the holomorphic
differential $dx/y$ on the canonical symplectic basis of homology
given by the ordering $\alpha_1,\alpha_2,\alpha_3$ of the roots of
$f$ (cf. \cite{mu}, Chapter IIIa, \S 5), and put $\tau :=
\omega_2/\omega_1$. Then we have the formulas
\[ \begin{array}{rcl}
\omega_1 \sqrt{ \alpha_1 - \alpha_3 } & = &  \pi \cdot
\vartheta(0;\tau)^2 \, ,  \\ \omega_1 \sqrt{ \alpha_1 - \alpha_2 }
& = &  \pi \cdot \vartheta(1/2;\tau)^2 \, ,   \\ \omega_1 \sqrt{
\alpha_2 - \alpha_3 } & = & \pi \cdot \exp(\pi i \tau /2)
    \cdot \vartheta(\tau/2;\tau)^2  \\
\end{array} \] for appropriate choices of the square roots. Let
$D:=16(\alpha_1-\alpha_2)^2(\alpha_1-\alpha_3)^2(\alpha_2-\alpha_3)^2
= p^3-27q^2$ be the discriminant of $f$. Then the formula
\[ D=(2\pi)^{12} \cdot \omega_1^{-12} \cdot \Delta(\tau) \] holds.
\end{lemma}
\begin{proof} The first set of formulas follows by an application
of Thomae's formula, \emph{cf.} \cite{mu}, Chapter IIIa, \S 5. The
other formula follows from the first and from the formula
\[ \left( \exp(\pi i \tau/4) \cdot \vartheta(0;\tau)
\vartheta(1/2;\tau) \vartheta(\tau/2;\tau) \right)^8 = 2^8 \cdot
\Delta(\tau) \] mentioned above.
\end{proof}
\begin{proposition} Let $X$ be a Riemann surface of genus 1 and suppose
that $y^2=4x^3-px-q=: f(x)$ is a Weierstrass equation for $X$.
Write $f(x) =4(x-\alpha_1)(x-\alpha_2)(x-\alpha_3)$. Let
$P_1=(\alpha_1,0)$, $P_2=(\alpha_2,0)$ and $P_3 =(\alpha_3,0)$.
Then the formulas
\[
G(P_1, P_2)^{12}  =  \frac{ 16 \cdot | \alpha_1 - \alpha_2 |^2 }
 { | \alpha_1 - \alpha_3| \cdot | \alpha_2 - \alpha_3 | } \, , \]
\[ G(P_1, P_3)^{12}  = \frac{ 16 \cdot | \alpha_1 - \alpha_3 |^2 }
 { | \alpha_1 - \alpha_2| \cdot | \alpha_3 - \alpha_2 | } \, , \]
\[ G(P_2, P_3)^{12} =  \frac{ 16 \cdot | \alpha_2 - \alpha_3 |^2 }
 { | \alpha_2 - \alpha_1| \cdot | \alpha_3 - \alpha_1 | }  \]
hold.
\end{proposition}
\begin{proof} This follows directly from Lemma \ref{classicfactsI} and
the explicit formula for $G(0,z)$ in Proposition \ref{explG}.
\end{proof}
We remark that this proposition has been obtained by Szpiro in
\cite{sz} in the special case that $X$ is the Riemann surface
associated to a Frey curve $y^2=x(x+a)(x-b)$, where $a,b$ are
non-zero integers with $2^4 | a$ and $b \equiv -1 \mod 4$
(\emph{cf.} \cite{sz}, Section 1.3).

\section{Arakelov projection formula} 

In this section we prove a
projection formula for Arakelov intersections on elliptic arithmetic
surfaces. The essential idea is to use 
the analytic projection formula from Theorem \ref{complexproj}; the
rest of the proof is quite straightforward. We will use the Arakelov
projection formula in Section \ref{average}.

Let $p :\ee
\to B=\mathrm{Spec}(O_K)$ be
an arithmetic surface over the ring of integers $O_K$ of a
number field $K$. Here and below we assume that $\ee$ is a regular
scheme. As in \cite{ar} we have on $\ee$
the notion of an Arakelov divisor: this is a formal sum of a Weil
divisor $D_\mathrm{fin}$ on $\ee$ and an infinite part
$D_\mathrm{inf}=\sum_\sigma \alpha_\sigma \cdot
E_\sigma$, the sum running over the complex embeddings of $K$, with
$\alpha_\sigma$ real numbers and with the $E_\sigma$ formal symbols
corresponding to the Riemann surfaces associated to the curves
$E \times_{K,\sigma} \cc$. The Arakelov divisors form a group
$\widehat{\mathrm{Div}}(\ee)$. 
To each non-zero rational function $f \in K(E)$ one
associates the corresponding Arakelov principal divisor $(f)$ with
$(f)_\mathrm{fin}$ the usual principal divisor associated to $f$, and
with $\alpha_\sigma$ given by $\alpha_\sigma = -\int_{E_\sigma} \log
|f|_\sigma \mu_\sigma$. Here $\mu_\sigma$ is the fundamental
(1,1)-form on $E_\sigma$. We denote
by $\widehat{\mathrm{Cl}}(\ee)$ the group of Arakelov divisors on
$\ee$ modulo the principal divisors. It was proved in \cite{ar} that 
there exists a natural bilinear symmetric intersection product
$(\cdot,\cdot)$ on the group of Arakelov divisors, factoring through
the principal divisors to give a natural bilinear symmetric
intersection pairing on $\widehat{\mathrm{Cl}}(\ee)$. The definition
of this intersection product is quite straightforward, except for the
crucial case of the intersection $(P,Q)$ of two sections $P,Q : B \to
\ee$ of $p$, which consists of a finite contribution
$(P,Q)_\mathrm{fin}$ given in the usual way, and an infinite
contribution $(P,Q)_\mathrm{inf}$ given as a sum $\sum_\sigma
(P,Q)_\sigma$ with $(P,Q)_\sigma = 
-\log G_\sigma(P_\sigma,Q_\sigma)$. Here
$G_\sigma$ the Arakelov-Green function on $E_\sigma$. 

Our Arakelov projection formula is a projection formula involving
pushforwards and pullbacks of Arakelov divisors, which we define as
follows.
\begin{definition} \label{pushfpullb} Let $p:\ee \to B$ and $p' : \ee'
\to B$ be elliptic arithmetic surfaces, and suppose there exists a
proper $B$-morphism $f : \ee \to \ee'$.
Let $D$ be an Arakelov
divisor on $\ee$, and write $D = D_{\mathrm{fin}} + \sum_{\sigma}
\alpha_{\sigma} \cdot E_{\sigma}$. The pushforward $f_*D$ of $D$
is defined to be 
the Arakelov divisor $f_*D := f_*D_{\mathrm{fin}} + d \cdot
\sum_{\sigma} \alpha_{\sigma} \cdot E'_{\sigma}$ on $\ee'$ where
$f_*D_{\mathrm{fin}}$ is the usual pushforward of the Weil divisor
$D_{\mathrm{fin}}$.
Next let $D'$ be an Arakelov divisor on
$\ee'$. The pullback $f^*D'$ of $D'$ is the Arakelov divisor
$f^*D' := f^*D'_{\mathrm{fin}} + \sum_{\sigma} \alpha'_{\sigma}
\cdot E_{\sigma}$ on $\ee$ 
where $f^*D'_{\mathrm{fin}}$ is the pullback of the Weil
divisor $D'_{\mathrm{fin}}$ on $\ee'$, defined in the usual way using
Cartier divisors.
\end{definition}
Our result is then as follows.
\begin{theorem} \label{projection} (Arakelov projection formula) Let
$E$ and $E'$ be elliptic curves defined over a number field $K$, and
let $p: \ee \to B$ and $p' : \ee' \to B$ be arithmetic
surfaces over the ring of integers of $K$ with generic fibers
isomorphic to $E$ and $E'$, respectively.
Suppose we have an isogeny
$f:E \to E'$, and suppose that $f$ extends to a $B$-morphism
$f:\ee \to \ee'$. Let $D$ be an Arakelov divisor on $\ee$ and let
$D'$ be an Arakelov divisor on $\ee'$. Then the equality of
intersection products $(f^*D',D)=(D',f_*D)$ holds.
\end{theorem}
\begin{proof} We may restrict to the case where both $D$ and $D'$ are Arakelov
divisors with trivial contributions ``at infinity''. By the moving
lemma on $\ee'$ (\emph{cf.} \cite{liu}, Corollary 9.1.10) 
we can find a function $g \in K(E')$ such
that $D'' := D' + (g)_{\mathrm{fin}}$ and $f_*D$ have no
components in common. Obviously $D''+(g)_{\mathrm{inf}}$ is
Arakelov linearly equivalent to $D'$, and hence by a computation as in
Theorem \ref{complexproj} the Arakelov divisor 
$f^*D''+(f^*g)_{\mathrm{inf}}$
is Arakelov linearly equivalent to $f^*D'$. It is therefore
sufficient to prove that $(f^*D'' + (f^*g)_{\mathrm{inf}},D) =
(D'' + (g)_{\mathrm{inf}}, f_*D)$. It is clear that $
((f^*g)_{\mathrm{inf}},D) = ((g)_{\mathrm{inf}},f_*D)$, so it
remains to prove that $(f^*D'',D)=(D'',f_*D)$. By the traditional
projection formula (\emph{cf.} \cite{liu}, Theorem 9.2.12 and Remark
9.2.13) we have $(f^*D'',D)_{\mathrm{fin}}
=(D'',f_*D)_{\mathrm{fin}}$. For the contributions at
infinity we can reduce to the case where $D$ and $D''$ are
sections of $\ee \to B$ and $\ee' \to B$, respectively. Let
$\sigma$ be a complex embedding of $K$. Let $D_{\sigma}$ and
$D''_{\sigma}$ be the points corresponding to $D$ and $D''$ on
$E_{\sigma}$ and $E'_{\sigma}$. Then for the local
intersection at $\sigma$ we have $(f^*D'',D)_\sigma =
(D'',f_*D)_\sigma$ by the analytic projection formula from
Proposition \ref{complexproj}. The theorem follows.
\end{proof}
\begin{remark} In general, an isogeny $f : E \to E'$ may not
extend to a morphism $\ee \to \ee'$. However, if $\ee'$ is a
minimal arithmetic surface (\emph{cf.} \cite{liu}, Section 9.3.2), 
then it contains the N\'eron model of
$E'/K$, and hence by the universal property of the N\'eron model,
any isogeny $f$ extends. In any case we can achieve that $f$
extends after blowing up finitely many closed points on $\ee$.
\end{remark}
The following corollary appears in Szpiro's
paper \cite{sz}.
\begin{corollary} (Szpiro \cite{sz})
Let $D_1,D_2$ be Arakelov divisors on $\ee'$.
Let $N$ be the degree of $f$. Then the formula \[ (f^*D_1, f^*D_2) = N \cdot
(D_1,D_2) \] holds.
\end{corollary}
\begin{proof} It is not difficult to see (\emph{cf.} \cite{liu},
Theorem 7.2.18 and Proposition 9.2.11)
that $f_*f^*D_2 = N \cdot D_2$. Theorem \ref{projection} then
gives $ (f^*D_1,f^*D_2) = (D_1,f_* f^* D_2) = (D_1, N \cdot D_2) =
N \cdot (D_1, D_2) $.
\end{proof}

\section{Self-intersection of a point}

Let $p:\ee \to B$ be an elliptic 
arithmetic surface. The image
of a section $P : B \to \ee$ gives rise to a divisor on $\ee$,
which we also denote by $P$. Given the framework of arithmetic
intersection theory, it is natural to ask for the
self-intersection $(P,P)$ of $P$. The question has been solved in
the case that $P$ is the zero section by Szpiro.
\begin{theorem} \label{(O,O)} (Szpiro \cite{sz})
Let $E$ be a semi-stable elliptic curve over a number field $K$,
and let $p: \ee \to B$ be its regular minimal model over the ring of
integers of $K$. Let $O : B \to \ee$ be the zero section of $p$,
and denote by $\Delta(E/K)$ the minimal discriminant ideal of
$E/K$. Then the formula
\[ (O,O) = -\frac{1}{12} \log |N_{K/\qq}(\Delta(E/K)| \] holds.
\end{theorem}
\begin{proof}
By Proposition \ref{adjunction} we need to compute the Arakelov degree 
$\Ardeg \, O^* \omega_{\ee/B}$, with $\omega_{\ee/B}$ the relative
dualising sheaf of $p: \ee \to B$. It is well-known that there
exists a canonical isomorphism $O^* \omega_{\ee/B} \isom p_*
\omega_{\ee/B}$. Now the line bundle $(p_*
\omega_{\ee/B})^{\otimes 12}$ contains a canonical global section
$\Lambda_{\ee/B}$ coming from the canonical isomorphism $(p_*
\omega)^{\otimes 12} \isom O(\Delta)$ on the moduli stack of
stable elliptic curves, with $\Delta$ the discriminant locus.
Considering then the canonical section $\Lambda_{\ee/B}$ in $(O^*
\omega_{\ee/B})^{\otimes 12}$ we compute its norm. First of all,
the finite places yield a contribution $\log
|N_{K/\qq}(\Delta(E/K))|$. Next consider a complex embedding
$\sigma$ of $K$. Suppose we have an identification $E_\sigma \cong
\cc/\zz+\tau_\sigma \zz$ with $\tau_\sigma$ in the complex upper
half plane. Let $y^2 =4x^3-g_{2\sigma}x-g_{3\sigma}$ be the
associated Weierstrass equation, where $x = \wp_\sigma(z)$ and $y
= \wp'_\sigma(z)$, with $\wp_\sigma$ the Weierstrass
$\wp$-function associated to the lattice $\zz+\tau_\sigma\zz$. We
then have $\Lambda_\sigma = D_\sigma \cdot (dx/y)^{\otimes 12}$
where $D_\sigma$ is the discriminant of the above Weierstrass
equation. Moreover $dx/y$ is identified with $dz$. We can now
compute $\|\Lambda_\sigma\|_{\mathrm{Ar}}$ as follows: first by
Lemma \ref{classicfactsI} we have $ D_\sigma = (2\pi)^{12} \cdot
\Delta(\tau_\sigma)$, and second by Proposition \ref{AX} we have
$\|dz\|_{\mathrm{Ar}}= \sqrt{\imtau} / ((2\pi) \cdot
\|\eta\|(E_\sigma)^2)$. We obtain that
$\|\Lambda_\sigma\|_{\mathrm{Ar}}=1$ and hence the infinite
contributions vanish. This gives the proposition.
\end{proof}
The proof given in \cite{sz} is much more involved. 
The above proof in fact answers a question raised in
\cite{sz} on the norm $\| \Lambda \|_{\mathrm{Ar}}$ of $\Lambda$
in $\Omega^{\otimes 12}$.

The following proposition shows that Theorem \ref{(O,O)} in fact
gives the general answer to our question.
\begin{proposition} \label{self} Let $E$ be an elliptic curve over
a number field $K$, and let $p : \ee \to B$ be the regular minimal model
of $E$ over the ring of integers of $K$. Let $O : B \to \ee$ be
the zero-section. Then for any section $P : B \to \ee$ of $\ee \to
B$ we have $(P,P)=(O,O)$.
\end{proposition}
For the proof, we make use of the following lemma.
\begin{lemma} \label{structureoom}
Let $E$ be an elliptic curve over a number field
$K$, and let $p : \ee \to B$ be the regular minimal model of $E$ over the
ring of integers of $K$. Let $\omega_{\ee/B}$ be the relative
dualising sheaf of $p: \ee \to B$. 
Then we can write $\omega_{\ee/B}=\sum_b
\lambda_b \ee_b + \sum_\sigma \alpha_\sigma F_\sigma$ as Arakelov
divisors on $\ee$, the first sum running over the closed points
$b$ of $B$, with $\ee_b$ denoting the fiber at $b$ and with
$\lambda_b$ certain rational numbers; the second sum runs over
the complex embeddings of $K$, with $\alpha_\sigma$ certain real
numbers.
\end{lemma}
\begin{proof} Since $
\omega_{\ee/B}$ restricts to the trivial sheaf on the generic
fiber there exists a vertical divisor $V$ on $\ee$ such that
$\omega_{\ee/B} \cong O_\ee(V)$ as invertible sheaves. 
Since $\ee$ is minimal, the
divisor $V$ is numerically effective (\emph{cf.} \cite{liu}, Corollary
9.3.26), which implies $(V,C) \geq 0$ for every irreducible
component $C$ of a closed fiber. But also by the adjunction
formula in the vertical fibers (\emph{cf.} \cite{liu}, Section
9.1.3) we have $(V,\ee_b)=2p_a(E)-2=0$ for each closed fiber $\ee_b$ of $\ee$, so in fact
$(V,C)=0$ for each $C$. Since the kernel of the intersection
product is generated by the multiples of the fibers, this implies
that $V= \sum_b \lambda_b \cdot \ee_b$, where the $\lambda_b$ are
certain rational numbers. The lemma follows immediately from this.
\end{proof}
\begin{proof}[Proof of Proposition \ref{self}] The adjunction formula
Proposition \ref{adjunction} shows that we need to prove that 
$\Ardeg \, P^*
\omega_{\ee/B}=\Ardeg \, O^*\omega_{\ee/B}$. But this is immediate 
from Lemma \ref{structureoom}.
\end{proof}
Note that Lemma \ref{structureoom} also proves that
$(\omega_{\ee/B},\omega_{\ee/B})=0$ on a minimal elliptic
arithmetic surface $p: \ee \to B$, a fact observed by Faltings in
\cite{fa} in the case of a semi-stable elliptic arithmetic
surface.

\section{Average height of quotients} \label{average}

In this final section we study the average height of quotients of
an elliptic curve by its cyclic subgroups of fixed order. Using
our results from the previous sections, we give an alternative
proof of a formula due to Autissier \cite{aut}. A slightly less
general result appears already in \cite{szpull}, and in fact our
method is very much in the spirit of this latter paper. The main
difference is perhaps that in our approach we do not need to
consider the distribution of torsion points on the bad fibers. 
In fact we do not need any non-trivial
arithmetic information at all; the main ingredients are the
Arakelov projection formula from Theorem \ref{projection}, the
formula for the ``energy of an isogeny'' from Theorem
\ref{kernel}, and the formula for the
self-intersection of a point from Theorem \ref{(O,O)}. 
Amusingly, we shall mention at the
end of this section how a purely arithmetic result, namely the
injectivity of torsion, follows from our Arakelov-theoretic
results.

We start with an explicit formula for $h_F(E)$. This formula is
certainly well-known,
\emph{cf.} \cite{sil}, Proposition 1.1.
\begin{proposition} \label{explicitheight} Let $E$ be a semi-stable
elliptic curve over a number field $K$. Let $\Delta(E/K)$
be the minimal discriminant ideal of $E/K$. Then the formula
\[  h_F(E) = \frac{1}{[K:\qq]} \left( \frac{1}{12} \log
|N_{K/\qq}(\Delta(E/K))| - \frac{1}{12} \sum_\sigma \log (
(2\pi)^{12} \|\Delta\|(E_\sigma)) \right) \] holds. Here the sum
runs over the complex embeddings of $K$.
\end{proposition}
\begin{proof} As is explained in the proof of Theorem \ref{(O,O)},
the line bundle $(p_* \omega_{\ee/B})^{\otimes 12}$ contains a
canonical section $\Lambda_{\ee/B}$, which has divisor given by
$\Delta(E/K)$ on $B$. This accounts for the finite contribution
$|N_{K/\qq}(\Delta(E/K))|$. Next, at a complex embedding $\sigma$
of $K$ we have $\Lambda_\sigma = D_\sigma \cdot (dx/y)^{\otimes
12}$ where $D_\sigma$ is the discriminant of a Weierstrass
equation $y^2 = 4x^3 -p_\sigma x -q_\sigma$ associated to
$E_\sigma$. Let $(\omega_{1\sigma}|\omega_{2\sigma})$ be a period
matrix of $dx/y$ on a canonical symplectic basis associated to an
ordering of the roots of $f$, and let $\tau_\sigma =
\omega_{2\sigma}/\omega_{1\sigma}$. By Lemma \ref{classicfactsI}
we have $ D_\sigma = (2\pi)^{12} \omega_{1\sigma}^{-12} \cdot
\Delta(\tau_\sigma)$, and by Riemann's second bilinear relations we
have $ \| dx/y \|_\sigma^2 = |\omega_{1\sigma}|^2 \cdot \imtau$.
Together this yields $\| \Lambda \|_\sigma = (2\pi)^{12} \cdot
\|\Delta\|(E_\sigma)$. This gives the infinite contribution to
$h_F(E)$.
\end{proof}
Now let's turn to the result of Autissier.
First we introduce some notations. Let $N$ be a positive integer. Then
we denote by $e_N$ the number of cyclic subgroups of order $N$ on an
elliptic curve defined over $\cc$, \emph{i.e.}
\[ e_N := N \prod_{p|N} \left( 1 + \frac{1}{p}
\right) \, , \] where the product is over the primes dividing $N$.
Further we put
\[ \lambda_N := \sum_{ {p|N \atop p^r \| N}} \frac{ p^r - 1 }{ p^{r-1} (p^2-1) } \log p
\, , \] where the notation $p^r \| N$ means that $p^r |N$ and
$p^{r+1} \nmid N$. For an elliptic curve $E$ and a finite subgroup $C$
of $E$ we denote by $E^C$ the quotient of $E$ by $C$.

In \cite{szpull} we find the following theorem.
\begin{theorem} (Szpiro-Ullmo, \cite{szpull}) Let $E$ be a semi-stable elliptic curve defined over a
number field $K$. Suppose that $E$ has no complex multiplication over
$\overline{K}$ and that the absolute Galois group
$\mathrm{Gal}(\overline{K}/K)$ acts transitively on the points of
order $N$ on $E$. Let $C$ be a cyclic subgroup of order $N$ on $E$.
Then the formula
\[ h_F(E^C) = h_F(E) + \frac{1}{2} \log N - \lambda_N \] holds.
\end{theorem}
One may wonder what one can say without the assumption that
$\mathrm{Gal}(\overline{K}/K)$ acts transitively. In \cite{aut} we
find a proof of the following statement. The price we pay for dropping
the assumption is that we
can only deal with the average over all $C$.
\begin{theorem} (Autissier \cite{aut}) \label{aut}
Let $E$ be an elliptic curve
defined over a number field $K$. Then the formula
\[ \frac{1}{e_N} \sum_C h_F(E^C) = h_F(E) + \frac{1}{2} \log N -
\lambda_N  \]  holds, where the sum runs over the cyclic
subgroups of $E$ of order $N$.
\end{theorem}
In fact, this formula was already stated in \cite{szpull} under
the restriction that $N$ is squarefree. Autissier's proof uses the
Hecke correspondence $T_N$ and a generalised intersection theory
for higher-dimensional arithmetic varieties. The disadvantage of
this approach is that the analytic machinery needed to deal with
the contributions at infinity becomes quite complicated. We will
give a proof of Theorem \ref{aut} which is much more elementary.
Besides this merit, we also think that the structure of the
somewhat strange constant $\lambda_N$ becomes more clear through
our approach. It would be interesting to have a generalisation of
Theorem \ref{aut} to abelian varieties of higher dimension.

Theorem \ref{aut} follows directly from the following two
propositions, by using the explicit formula for $h_F$ in
Proposition \ref{explicitheight}. The next proposition occurs as Lemme
5.4 in \cite{szpull}.
\begin{proposition} \label{firstprop}
Let $E$ be a semi-stable elliptic curve over a number field $K$
and suppose that all $N$-torsion points are $K$-rational. Then one
has
\[ \sum_C \left(  \log
|N_{K/\qq}(\Delta(E/K))| - \log
|N_{K/\qq}(\Delta(E^C/K)| \right) = 0 \, . \] Here the sum runs
over the cyclic subgroups of $E$ of order $N$.
\end{proposition}
\begin{proposition} \label{secondprop}
Let $X$ be a Riemann surface of genus 1. Then
\[ \frac{1}{e_N} \sum_C \left( \frac{1}{12} \log \|\Delta\|(X) -
\frac{1}{12} \log \|\Delta\|(X^C) \right) = \frac{1}{2} \log N -
\lambda_N \, ,
\] where the sum runs over the cyclic subgroups of $X$ of order $N$.
\end{proposition}
Our first step is to
reduce these two propositions to the following two:
\begin{proposition} \label{main1}
Let $E$ be a semi-stable elliptic curve over a number field $K$
and suppose that all $N$-torsion points are $K$-rational. Extend
all $N$-torsion points of $E$ over the regular minimal model of
$E/K$. Then one has
\[ \sum_C \sum_{ {Q \in C
\atop Q \neq O} } (Q,O) = 0 \, , \] where the first sum runs over
the cyclic subgroups of $E$ of order $N$, and the second sum runs
over the non-zero points in $C$.
\end{proposition}
\begin{proposition} \label{main2}
Let $X$ be a Riemann surface of genus 1. Then one
has \[ \frac{1}{e_N} \sum_C \sum_{ {Q \in C \atop Q \neq 0} } \log
G(Q,0) = \lambda_N \, . \] Here the first sum runs over the cyclic
subgroups of $X$ of order $N$, and the second sum runs over the
non-zero points in $C$.
\end{proposition}
The latter proposition is an improvement of Proposition 6.5 in
\cite{szpull}, which gives an analogous statement, but only with
the left hand side summed over the complex embeddings of $K$, and
divided by $[K:\qq]$. Our result holds in full generality for an
arbitrary Riemann surface of genus 1.
\begin{proof}[Proof of Proposition \ref{firstprop} from Proposition \ref{main1}]
Let $C$ be any cyclic subgroup of $E$, and let $O'$ be the
zero-section of $E^C$. Extend it over the minimal regular model of
$E^C/K$. We then have 
\[ \frac{1}{12} \log
|N_{K/\qq}(\Delta(E/K))| - \frac{1}{12} \log
|N_{K/\qq}(\Delta(E^C/K)| = (O',O') - (O,O)   \] 
by Theorem \ref{(O,O)}. The latter is
equal to $ \sum_{ {Q \in C
\atop Q \neq O} } (Q,O) $  by
Theorem \ref{projection}. Summing over all cyclic
subgroups of $E$ of order $N$ and using
Proposition \ref{main1} we find the result.
\end{proof}
\begin{proof}[Proof of Proposition \ref{secondprop} from Proposition \ref{main2}]
By Theorem \ref{kernel} we have for any subgroup
$C$ of $X$ of order $N$ that
\[ \frac{1}{12} \log \|\Delta\|(X) - \frac{1}{12} \log
\|\Delta\|(X^C) = \frac{1}{2} \log N - \sum_{ {Q \in C \atop Q
\neq 0} } \log G(Q,0) \, .
\] The statement of Proposition \ref{secondprop} is then
immediate from Proposition \ref{main2}.
\end{proof}
In order to prove Proposition \ref{main1}, we make use of the
following combinatorial lemma.
\begin{lemma} \label{cyclic}
Let $M$ be a positive integer with $M | N$. Let $E$ be an elliptic
curve defined over an algebraically closed field of characteristic
zero. Then each cyclic subgroup of $E$ of order $M$ is contained
in exactly $e_N/e_M$ cyclic subgroups of order $N$.
\end{lemma}
\begin{proof} This follows easily by fixing a basis for the $N$-torsion and then
considering the induced natural transitive left action of
$\mathrm{SL}(2,\zz)$ on the set of cyclic subgroups of order $M$
and of order $N$.
\end{proof}
We may argue then as follows.
\begin{proof}[Proof of Proposition \ref{main1}] Let
$\overline{E}[M]$ be the set of points of exact order $M$ on $E$.
By Lemma \ref{cyclic} we have
\[ \sum_C \sum_{ {Q \in C
\atop Q \neq O} } (Q,O) = \sum_{ { M|N \atop M > 1}}
\frac{e_N}{e_M} \sum_{ Q \in \overline{E}[M]} (Q,O) \, .
\] We claim that for any positive integer $M$, we have $\sum_{Q
\in \overline{E}[M]}(Q,O) =0$. Indeed, we have \[ \sum_{Q \in
E[M], Q \neq O} (Q,O) = 0 \] for all $M$ by Theorem
\ref{projection} and then the claim follows by M\"obius inversion.
\end{proof}
Also for the proof of Proposition \ref{main2} we will need a
lemma. For a Riemann surface $X$ of genus 1, and $M>1$ an integer,
we put
\[ t(M) := \sum_{Q \in \overline{X}[M]}
\log G(Q,0) \, , \] the sum running over the set $\overline{X}[M]$
of points of exact order $M$ on $X$.

Part of the following lemma is also given in \cite{szpull}, \emph{cf.}
Lemme 6.2.
\begin{lemma} \label{ont}
We have \[ t(p^r)=\log p \] for any prime integer $p$ and any
positive integer $r$. Moreover we have $t(M)=0$ for any positive
integer $M$ which is not a prime power.
\end{lemma}
\begin{proof} By Corollary \ref{Ntorsion} we have
\[ \sum_{Q \in X[M], Q \neq 0} \log G(Q,0) = \log M \, . \]
The lemma follows from this by M\"obius inversion.
\end{proof}
\begin{proof}[Proof of Proposition \ref{main2}] For any divisor $M|N$,
let $\overline{X}[M]$ be the set of points of exact order $M$ on
$X$ and let $t(M) = \sum_{Q \in \overline{X}[M]}
\log G(Q,0)$ as in Lemma \ref{ont} where it is understood that
$t(1) = 0$. Then by Lemma
\ref{cyclic} we can write
\[ \frac{1}{e_N} \sum_C \sum_{ {Q \in C \atop Q \neq 0} } \log
G(Q,0) = \frac{1}{e_N} \sum_{M|N} \frac{e_N}{e_M} \cdot t(M) \, .
\] Lemma \ref{ont} gives us that
\[  \frac{1}{e_N} \sum_{M|N} \frac{e_N}{e_M} \cdot t(M) =
\sum_{ { p|N \atop p^r \| N } } \left( \frac{1}{e_p} + \cdots +
\frac{1}{e_{p^r}} \right) \log p \, . \] Finally note that
$e_{p^k} = p^k(1+1/p)$ which gives
\[ \frac{1}{e_p} + \cdots +
\frac{1}{e_{p^r}} = \frac{ p^r -1 }{ p^{r-1} (p^2-1) } \, . \]
From this the result follows.
\end{proof}
\begin{remark} An alternative proof of Proposition
\ref{secondprop} can be given by classical methods using modular
forms identities, see for instance \cite{cata}, Proposition
VII.3.5(b) for the case that $N$ is a prime, and \cite{aut}, Lemme
2.2 and Lemme 2.3 for the general case. We preferred to give an
argument using Arakelov theory, indicating that Arakelov theory
can sometimes be used to derive analytic results on Riemann
surfaces in a short and clean manner.
\end{remark}
We finish with a corollary from the results above. This corollary
gives another interpretation to the constant $\lambda_N$.
\begin{corollary} \label{finitelambda}
Let $E$ be a semi-stable elliptic curve over a number field $K$
and suppose that all $N$-torsion points are $K$-rational. Extend
these torsion points over the minimal regular model of $E/K$. Then
one has
\[ \frac{1}{[K:\qq]} \frac{1}{e_N} \sum_C \sum_{ {Q \in C \atop Q
\neq O} } (Q,O)_{\mathrm{fin}} = \lambda_N \, ,
\] where the first sum runs over the cyclic subgroups of $E$ of
order $N$, and the second sum runs over the non-zero points in
$C$.
\end{corollary}
\begin{proof}
Let $C$ be a finite cyclic subgroup of $E$. Note that by definition of the
Arakelov intersection product
\[ \sum_{ {Q \in C
\atop Q \neq O} } (Q,O) = \sum_{ {Q \in C \atop Q \neq O} }
(Q,O)_{\mathrm{fin}} -  \sum_{ {Q \in C \atop Q \neq O} }
\sum_\sigma \log G(Q^\sigma,0) \, . \] The corollary follows
therefore
easily from Proposition \ref{main1} and Proposition \ref{main2}.
\end{proof}
Note that Corollary \ref{finitelambda} is purely arithmetical in
nature. It should also be possible to give a direct proof, but
probably this would require a more \emph{ad hoc} approach, making
for instance a case distinction between the supersingular and the
non-supersingular primes for $E/K$. Also note that Corollary
\ref{finitelambda} immediately gives the classical arithmetic
result that, for any prime number $p$, the $p$-torsion points are injective on a fiber at a
prime of characteristic different from $p$. Indeed, 
take $N=p$ in the
formula from Corollary \ref{finitelambda}, then the right hand
side is a rational multiple of $\log p$, and so the same holds for
the left hand side. This means that the local intersections
$(Q,O)_{\mathrm{fin}}$, which are always non-negative, 
are in fact zero at primes of characteristic
different from $p$. Hence, each $p$-torsion point
$Q$ stays away from $O$ on
fibers above such primes. Of course the
argument can be repeated with $O$ replaced by any other $p$-torsion
point.

\subsection*{Acknowledgments} The author wishes to thank his thesis
advisor Gerard van der Geer for his encouragement and helpful
remarks. Also he thanks Professor Qing Liu for pointing out a mistake
in an earlier version of this note.

\end{document}